\documentclass[12pt]{amsart}
\usepackage{amsfonts,amsmath,amsthm,amscd,amssymb,latexsym,cite,verbatim,texdraw,floatflt,
caption2}
\RequirePackage{hyperref}

\usepackage[a4paper]{geometry}

\theoremstyle
{plain}

\keywords{coarse structure,  ballean, vector balleans, free vector ballean.}
\date{\today}

\address{I. Protasov, Department of Computer Science and Cybernetics, Kyiv University, Volodymyrska 64, 01033, Kyiv, Ukraine}
\email{i.v.protasov@gmail.com; }
\address{K. Protasova, Department of Computer Science and Cybernetics, Kyiv University, Volodymyrska 64, 01033, Kyiv, Ukraine}
\email{ksuha@freenet.com.ua}

\begin{document}

\title{A note on free vector balleans}

\author{Igor Protasov and Ksenia Protasova}

\maketitle
\vskip 5pt

{\bf Abstract.}
A vector balleans is a  vector  space over $\mathbb{R}$  endowed with a coarse structure in such a way that  the vector  operations are coarse mappings. We prove that, for every ballean  $(X, \mathcal{E})$, there exists the unique free vector ballean $\mathbb{V}(X, \mathcal{E})$ and  describe the coarse structure of $\mathbb{V}(X, \mathcal{E})$.
It is shown that normality of $\mathbb{V}(X, \mathcal{E})$ is equivalent to metrizability  of
$(X, \mathcal{E})$.
\vskip 10pt

{\bf MSC: } 46A17, 54E35.
\vskip 5pt

{\bf Keywords:} coarse structure,  ballean, vector ballean, free vector ballean.

\section{Introduction}

Let $X$  be a set. A family $\mathcal{E}$ of subsets of $X\times X$ is called a {\it coarse structure } if

\begin{itemize}
\item{}   each $\varepsilon\in \mathcal{E}$  contains the diagonal  $\bigtriangleup _{X}$,
$\bigtriangleup _{X}= \{(x,x): x\in X\}$;

\item{}  if  $\varepsilon,  \delta\in \mathcal{E}$ then $\varepsilon \circ\delta\in\mathcal{E}$ and
$\varepsilon^{-1}\in \mathcal{E}$,   where    $\varepsilon \circ\delta=\{(x,y): \exists z((x,z) \in\varepsilon,  \   \ (z, y)\in\delta)\}$,   $\varepsilon^{-1}=\{(y,x): (x,y)\in\varepsilon\}$;

\item{} if $\varepsilon\in\mathcal{E}$ and $\bigtriangleup_{X}\subseteq \varepsilon^{\prime}\subseteq\varepsilon $   then
$\varepsilon^{\prime}\in \mathcal{E}$;

\item{}  for any   $x,y\in X$, there exists $\varepsilon\in \mathcal{E}$   such that $(x,y)\in\varepsilon$.

\end{itemize}
\vskip 5pt

A subset $\mathcal{E}^{\prime} \subseteq \mathcal{E}$  is called a
{\it base} for $\mathcal{E}$  if, for every $\varepsilon\in \mathcal{E}$, there exists
  $\varepsilon^{\prime}\in \mathcal{E}^{\prime}$  such  that
  $\varepsilon\subseteq\varepsilon^{\prime}$.
For $x\in X$,  $A\subseteq  X$  and
$\varepsilon\in \mathcal{E}$, we denote
$\varepsilon[x]= \{y\in X: (x,y) \in\varepsilon\}$,
 $\varepsilon[A] = \cup_{a\in A}   \   \   \varepsilon[a]$
 and say that  $\varepsilon[x]$
  and $\varepsilon[A]$
   are {\it balls of radius} $\varepsilon$
   around $x$  and $A$.

The pair $(X,\mathcal{E})$ is called a {\it coarse space} \cite{b11}  or a ballean \cite{b8}, \cite{b10}.

Each subset $Y\subseteq X$  defines the {\it subballean}
$(Y, \mathcal{E}_{Y})$, where $\mathcal{E}_{Y}$  is the restriction of $\mathcal{E}$  to $Y\times Y$.
A subset $Y$  is called {\it bounded} if $Y\subseteq \varepsilon[x]$  for some $x\in X$  and
$\varepsilon\in \mathcal{E}$.

Let $(X, \mathcal{E})$, $(X^{\prime}, \mathcal{E}^{\prime})$   be balleans.
A mapping $f: X\longrightarrow  X^{\prime}$  is called
{\it coarse } or {\it macro-uniform} if,  for every  $\varepsilon\in\mathcal{E}$,
there  exists $\varepsilon^{\prime}\in\mathcal{E}^{\prime}$
 such that $f(\varepsilon[x])\subseteq \varepsilon^{\prime}[f(x)]$
  for each $x\in X$.
If $f$ is a bijection such that $f$  and $f^{-1}$  are coarse then $f$  is called an {\it asymorphism}.

Every metric $d$ on a set $X$ defines the {\it metric ballean }  $(X, \mathcal{E} _{d})$,
 where $\mathcal{E} _{d}$  has the base $\{\{(x,y): d(x,y)< r  \} :  r\in\mathbb{R}^{+}\}$.
We say  that a ballean $(X, \mathcal{E})$  is {\it metrizable}  if there exists a
 metric $d$  on $X $  such that  $\mathcal{E}= \mathcal{E} _{d}$.
In what follows, we consider $\mathbb{R}$  as a ballean defined by the metric $d(x, y)= |x-y|$.

Given two balleans
$(X_{1}, \mathcal{E}_{1})$,  $(X_{2}, \mathcal{E}_{2})$,
we defines the product $(X_{1} \times  X_{2} ,  \mathcal{E} )$,  where $\mathcal{E}$  has the base $\mathcal{E}_{1}\times\mathcal{E}_{2}$.

Let $\mathcal{V}$ be a vector  space over $\mathbb{R}$  and let $\mathcal{E}$ be a coarse structure on $\mathcal{V}$.
Following \cite{b5}, we say that  $(\mathcal{V}, \mathcal{E})$  is a {\it vector ballean} if the operation
$$\mathcal{V} \times \mathcal{V} \longrightarrow  \mathcal{V},  \  \  (x,y) \longmapsto   x+y  \   \   and   \   \   \mathbb{R}\times \mathcal{V},  \  \   (\lambda, x)\longmapsto \lambda x $$   are coarse.

A family  $\mathcal{I}$  of subsets of $\mathcal{V}$  is called a {\it vector ideal} if

$(1) \  \  $   if   $U, U^{\prime} \in \mathcal{I}$  and
$W\subseteq  U$  then $U\cup U^{\prime} \in \mathcal{I}$  and $W\in \mathcal{I}$;

$(2) \  \  $    for every  $U\in \mathcal{I}$,   $U +U\in \mathcal{I}$;

$(3) \  \  $   for any  $U\in \mathcal{I}$  and $\lambda\in\mathbb{R}^{+}$,  $[-\lambda, \lambda]  \ U\in\mathcal{I}$,
 where  $[-\lambda, \lambda] \ U = \cup \  \{\lambda^{\prime} \ U : \lambda^{\prime}\in [-\lambda, \lambda]\}$;

$(4) \  \  $   $\cup \ \mathcal{I} = \mathcal{V}$.

A family $\mathcal{I}^{\prime}\subseteq \mathcal{I}$
is called a {\it base}  for $\mathcal{I}$  if, for each  $U\in \mathcal{I}$,   there  is  $U^{\prime}\in\mathcal{I}^{\prime}$
such that
$U \subseteq U^{\prime}$.

If $ (\mathcal{V}, \mathcal{E})$ is a vector  ballean then the family $\mathcal{I}$
of all bounded subsets of $\mathcal{V}$ is a vector ideal. On the
 other hand, every vector ideal $\mathcal{I}$  on $\mathcal{V}$  defines the
 vector ballean  $(\mathcal{V}, \mathcal{E})$, where $\mathcal{E}$ is a coarse structure
 with the base
 $\{\{(x, y): x-y\in U\}:  U\in\mathcal{I}\}$.
Thus, we have got a bijective correspondence between vector balleans  on $\mathcal{V}$ and vector ideas.
Following this correspondence, we  write
$(\mathcal{V}, \mathcal{I})$ in place of  $ (\mathcal{V}, \mathcal{E})$.

Let
$(\mathcal{V}, \mathcal{I})$, $(\mathcal{V}^{\prime}, \mathcal{I}^{\prime})$
  be vector  balleans. We note that  a linear mapping
  $f: \mathcal{V}\longrightarrow \mathcal{V}^{\prime}$  is coarse if and only if
   $f(U)\in\mathcal{I}^{\prime}$
    for each
    $U\in\mathcal{I}$.

In Section 2, we show that, for every baleen  $(X, \mathcal{E})$,  there exists the unique vector ideal
$\mathcal{I}_{(X,\mathcal{E})}$ on the
 vector space $\mathbb{V}(X)$ with the basis $X$ such  that

  \begin{itemize}
\item{}   $(X, \mathcal{E})$ is a subballean of  $(\mathbb{V} (X),  \mathcal{I}_{(X,\mathcal{E})})$;

\item{}    for every  vector ballean $(\mathcal{V}, \mathcal{I})$,  every coarse
 mapping
 $(X, \mathcal{E})\longrightarrow (\mathcal{V}, \mathcal{I})$
   gives rise to the unique coarse linear mapping
    $(\mathbb{V} (X),  \  \mathcal{I}_{(X,\mathcal{E})})\longrightarrow (\mathcal{V}, \mathcal{I})$ .
 \vskip 5pt

\end{itemize}

We denote $\mathbb{V} (X,  \mathcal{E})= (\mathbb{V} (X),  \mathcal{I}_{(X,\mathcal{E})})$
 and say that $\mathcal{I}_{(X,\mathcal{E})}$
  and $\mathbb{V} (X,  \mathcal{E})$
  are {\it free vector ideal} and {\it free vector ballean } over $( X, \mathcal{E})$.

Free vector balleans can be considered as  counterparts of  free vector spaces studied in many papers, for examples, \cite{b1}, \cite{b3},  \cite{b4}.
It should be mentioned that the
 free activity  in topological algebra was initiated by the  famous paper of Markov  on free  topological groups \cite{b6}.
For free coarse groups see \cite{b9}.

\section{Construction}

Given a ballean $( X, \mathcal{E})$, we consider $X$  as the  basis of the vector space
$\mathbb{V} (X)$.
For each $\varepsilon\in\mathcal{E}$
 and  $n\in\mathbb{N}$,
 we set
 $\mathcal{D}_{\varepsilon}= \{x-y: (x, y)\in\varepsilon\}$
  and denote by
  $\mathcal{S}_{n,\varepsilon}$
    the sum of  $n$ copies of
    $[-n, n]\mathcal{D}_{\varepsilon}$.\vskip 5pt

{\bf Theorem 1.  } {\it
Let $( X, \mathcal{E})$ be a ballean and let  $z\in X$.
Then  the family
$$\{\mathcal{S}_{n,\varepsilon} +  [-n, n]z: \varepsilon\in\mathcal{E}, n\in\mathbb{N} \}$$
  is a base of the free vector ideal
$\mathcal{I}_{(X,\mathcal{E})}$.
\vskip 5pt

Proof.}
We denote by  $\mathcal{I}$  the family of all  subsets
$U$  of $\mathbb{V} (X)$
 such that $U$  is  contained in some
 $\mathcal{S}_{n,\varepsilon} +  [-n, n]z$.
Clearly,  $\mathcal{I}$  satisfies  (1),  (2), (3)  from the definition of a vector ideal.
To see that
$\cup\mathcal{I} = \mathbb{V} (X)$,
 we take an arbitrary  $y\in X$   and choose  $\varepsilon\in \mathcal{E}$  so that
 $(y,z)\in \varepsilon$.
Then
$y=(y-z) + z \in \mathcal{D}_{\mathcal{E}} + z$.
In view of (2), (3),  we conclude that  $\cup \mathcal{I} = \mathbb{V} (X)$.

To  show that $(X,\mathcal{E})$ is a subballean of
$(\mathbb{V} (X), \ \mathcal{I}) $,
 we  denote by  $\mathcal{E}^{\prime}$  the coarse structure of the ballean
 $(\mathbb{V} (X), \ \mathcal{I} )$.
Since  $\mathcal{D}_{\mathcal{E}}\in\mathcal{I}$
for each $\varepsilon\in\mathcal{E}$, $\mathcal{E}= \mathcal{E}^{-1}$
 we have $\mathcal{E} \subseteq \mathcal{E}^{\prime}|_{X}$.
To verify the inclusion
$ \mathcal{E}^{\prime} \mid_{X}\subseteq \mathcal{E}$
 we take  $x,y \in  X$, assume that   $x-y\in \mathcal{S} _{n,\mathcal{E}}  +  [-n, n] z$ and show that  $(x, y)\in   \varepsilon ^{n}$.

 We write $x-y  =  \lambda_{1}  (x_{1}- y_{1}) +\ldots  + \lambda_{n}  (x_{n}- y_{n}) + \lambda_{n+1} z,
  (x_{i}, y_{i}) \in\varepsilon $,  $\lambda_{i}\in [-n, n]$.
Since
$( \lambda_{1} - \lambda_{1}) + \ldots  + (\lambda_{n}- \lambda_{n})  + \lambda_{n+1} = 0$,
we have $\lambda_{n+1}=0$ so
$x-y  =  \lambda_{1}  (x_{1}- y_{1}) +\ldots  + \lambda_{n}  (x_{n}- y_{n}) $.
If $(x_{1}, y_{1}), \ldots  , x_{n}, y_{n})\in\{x,y\}$
 then the statement is evident because
 $(x_{i},y_{i})\in \varepsilon$.
Assume that there  exists
$a\in \{x_{1}, y_{1}, \ldots  , x_{n}, y_{n}\}$
 such that
 $a\notin\{x,y\}$.
We take all items
$\lambda_{i} (x_{i}- y_{i})$, $i\in I$
 such that  $a\in \{x_{i}, y_{i}\}$,
  and denote by  $s$  the sum of all  these items.
The coefficient before $a$ in the canonical decomposition of $s$ by  the basis $X$  must be $0$.
We take
$\lambda_{k} (x_{k}- y_{k})$, $k\in I$.
If $x_{k}=a$
 or $y_{k}=a$
 then we replace each $a$  in
 $\lambda_{i} (x_{i}- y_{i})$, $i\in I$,
 to $y_{k}$ or $x_{k}$ respectively.
Then we get
$x-y  =  \lambda_{1}  (x_{1}^{\prime}- y_{1}^{\prime}) +\ldots  + \lambda_{n}  (x_{n}^{\prime}- y_{n}^{\prime}) $,
$a\notin \{x_{1}^{\prime}, y_{1}^{\prime}, \ldots  ,  x_{n}^{\prime}- y_{n}^{\prime}\}$
 and $(x_{i}^{\prime}, y_{i}^{\prime})\in \varepsilon^{2}$.
Repeating this trick, we run into the case
$x_{1}, y_{1}, \ldots  , x_{n}, y_{n}\in\{x,y\}$ and $(x_{i}, y_{i})\in \varepsilon^{n}$.

To conclude the proof, we observe that $\mathcal{I}$  is the minimal vector ideal on
$\mathbb{V}(X)$  such that $(X,\mathcal{E})$ is a subballean  of
$(\mathbb{V}(X), \mathcal{I})$.
If $(\mathcal{V}, \mathcal{I}^{\prime})$ is a ballean and
$f: (X, \mathcal{E})\longrightarrow (\mathcal{V}, \mathcal{I}^{\prime})$
 is a coarse mapping then the linear extension
 $h: (\mathbb{V}(X), \mathcal{I}) \longrightarrow $  $(\mathcal{V}, \mathcal{I}^{\prime})$
   of $f$  is coarse because  $h^{-1}(\mathcal{I}^{\prime})$  is a vector ideal on  $\mathbb{V}(X)$
   and  $\mathcal{I}\subseteq h^{-1}(\mathcal{I}^{\prime})$.
    $ \ \ \ \Box$

\section{Metrizability and normality}

{\bf Theorem 2.  } {\it
A ballean $\mathbb{V}(X, \mathcal{E})$
  is metrizable if and only if $(X, \mathcal{E})$ is  metrizable.
\vskip 5pt

Proof.} By [10,  Theorem 2.1.1], $(X, \mathcal{E})$  is  metrizable  if  and only if
 $\mathcal{E}$  has a  countable  base.
Apply Theorem 1.   $ \ \ \ \Box$ \vskip 5pt

Let $(X, \mathcal{E})$ be a ballean,  $A \subseteq X$.
A subset  $U$ of $X$ is called an {\it asymptotic neighbourhood}  of $A$  if, for every
$\varepsilon\in \mathcal{E}$, the set  $\varepsilon[A] \backslash U$  is bounded.

Two subset $A, B$  of $X$  are called {\it asymptotically disjoint (asymptotically  separated)}
if, for every $\varepsilon\in\mathcal{E}$,
the   intersection  $\varepsilon[A]\cap   \varepsilon[B]$  is bounded ($A, B$
 have disjoint  asymptotic  neighbourhoods).

A ballean $(X, \mathcal{E})$  is called {\it normal} \cite{b7} if any  two  asymptotically  disjoint  subsets of  $X$  are asymptotically  separated.
Every metrizable ballean is normal.

Given an arbitrary ballean  $(X, \mathcal{E})$,  the family $\mathcal{B}_{X}$  of all bounded subsets of
 $X$ is called a {\it bornology}  of $(X, \mathcal{E})$.
A subfamily $\mathcal{B}^{\prime}\subseteq  \mathcal{B}_{X}$  is called a {\it base}  of
  $\mathcal{B}_{X}$  if, for every
  $B\in\mathcal{B}_{X}$,  there exists $B^{\prime}\in\mathcal{B}_{X}$ such that $B\subseteq B^{\prime}$.
The minimal cardinality of bases  of
$\mathcal{B}_{X}$ is denoted by $cof \ \mathcal{B}_{X}$.

\vskip 7pt

{\bf Theorem 3.  } {\it
For every ballean $(X, \mathcal{E})$,
the free vector ballean
$\mathbb{V}(X, \mathcal{E})$
 is normal if and only if $\mathbb{V}(X, \mathcal{E})$ is metrizable.
\vskip 7pt

Proof.} For $|X|=1$, the statement is evident.
Let $|X|> 1$, $a\in X$,  $L = \mathbb{R}a$,  $Y= X \setminus\{a\}$.
Applying Theorem 1,  we conclude
 that the canonical isomorphism between
 $\mathbb{V}(X, \mathcal{E})$
  and
  $L \times \mathbb{V}(Y, \mathcal{E}_{Y})$
  is an asymorphism. If
  $\mathbb{V}(X, \mathcal{E})$
   is normal  then, by Theorem 1.4 from [2],
   $cof \ \mathcal{B}_{L}= cof \ \mathcal{B}_{Z}$,
     where $Z= \mathbb{V}(Y, \mathcal{E}_{Y})$.
Since $cof \ \mathcal{B}_{L}= \aleph_{0}$, $\mathcal{B}_{Z}$
has a countable base. To conclude the proof,
it suffices to note that
$\mathcal{B}_{Z}$
is the vector ideal such that
$\mathbb{V}(Y, \mathcal{E}_{Y})= (\mathbb{V} (Y), \mathcal{I})$.
Hence $\mathcal{I}$ has a countable base and
$\mathbb{V}(Y, \mathcal{E}_{Y})$
 is metrizable by Theorem 2.1.1. from [10].
 $ \ \  \  \Box$

\end{document}